\documentclass{amsart}
\usepackage{graphics}
\usepackage{amsmath}
\usepackage{latexsym}
\usepackage{amssymb}
\usepackage[all]{xy}
\xyoption{matrix}
\xyoption{arrow}

\begin{document}

\renewcommand{\th}{\operatorname{th}\nolimits}
\newcommand{\rej}{\operatorname{rej}\nolimits}
\newcommand{\extto}{\xrightarrow}
\renewcommand{\mod}{\operatorname{mod}\nolimits}
\newcommand{\ul}{\underline}
\newcommand{\Sub}{\operatorname{Sub}\nolimits}
\newcommand{\ind}{\operatorname{ind}\nolimits}
\newcommand{\exc}{\operatorname{exc\, }\nolimits}
\newcommand{\Fac}{\operatorname{Fac}\nolimits}
\newcommand{\add}{\operatorname{add}\nolimits}
\newcommand{\soc}{\operatorname{soc}\nolimits}
\newcommand{\Hom}{\operatorname{Hom}\nolimits}
\newcommand{\Rad}{\operatorname{Rad}\nolimits}
\newcommand{\RHom}{\operatorname{RHom}\nolimits}
\newcommand{\uHom}{\operatorname{\underline{Hom}}\nolimits}
\newcommand{\End}{\operatorname{End}\nolimits}
\renewcommand{\Im}{\operatorname{Im}\nolimits}
\newcommand{\Ker}{\operatorname{Ker}\nolimits}
\newcommand{\Coker}{\operatorname{Coker}\nolimits}
\newcommand{\Ext}{\operatorname{Ext}\nolimits}
\newcommand{\op}{{\operatorname{op}}}
\newcommand{\Ab}{\operatorname{Ab}\nolimits}
\newcommand{\id}{\operatorname{id}\nolimits}
\newcommand{\pd}{\operatorname{pd}\nolimits}
\newcommand{\A}{\operatorname{\mathcal A}\nolimits}
\newcommand{\C}{\operatorname{\mathcal C}\nolimits}
\newcommand{\D}{\operatorname{\mathcal D}\nolimits}
\newcommand{\X}{\operatorname{\mathcal X}\nolimits}
\newcommand{\Y}{\operatorname{\mathcal Y}\nolimits}
\newcommand{\F}{\operatorname{\mathcal F}\nolimits}
\newcommand{\Z}{\operatorname{\mathbb Z}\nolimits}
\renewcommand{\P}{\operatorname{\mathcal P}\nolimits}
\newcommand{\T}{\operatorname{\mathcal T}\nolimits}
\newcommand{\G}{\Gamma}
\renewcommand{\L}{\Lambda}
\newcommand{\bdot}{\scriptscriptstyle\bullet}
\renewcommand{\r}{\operatorname{\underline{r}}\nolimits}
\newtheorem{lem}{Lemma}[section]
\newtheorem{prop}[lem]{Proposition}
\newtheorem{cor}[lem]{Corollary}
\newtheorem{thm}[lem]{Theorem}
\newtheorem*{thmA}{Theorem A.1}
\newtheorem*{cor1}{Corollary A.2}
\newtheorem*{cor2}{Corollary A.3}
\newtheorem{rem}[lem]{Remark}
\newtheorem{defin}[lem]{Definition}

\title[Clusters and seeds]{Clusters and seeds in acyclic cluster algebras \\ \vspace{8mm}
\small{with an appendix coauthored in addition by \\ P. Caldero and B. Keller}}

\author[Buan]{Aslak Bakke Buan}
\address{Institutt for matematiske fag\\
Norges teknisk-naturvitenskapelige universitet\\
N-7491 Trondheim\\
Norway}
\email{aslakb@math.ntnu.no}

\author[Marsh]{Robert J. Marsh}
\address{Department of Mathematics \\
University of Leicester \\
University Road \\
Leicester LE1 7RH \\
England
}
\email{rjm25@mcs.le.ac.uk}

\author[Reiten]{Idun Reiten}
\address{Institutt for matematiske fag\\
Norges teknisk-naturvitenskapelige universitet\\
N-7491 Trondheim\\
Norway}
\email{idunr@math.ntnu.no}

\author[Todorov]{Gordana Todorov}
\address{Northeastern University\\
Department of Mathematics\\
360 Huntington Avenue\\
Boston, MA 02115\\
USA
}
\email{todorov@neu.edu}



\begin{abstract}
Cluster algebras are commutative algebras that were
introduced by Fomin and Zelevinsky in order to model the dual canonical
basis of a quantum group and total positivity in algebraic groups. Cluster
categories were introduced as a representation-theoretic model for
cluster algebras. In this article we use this representation-theoretic
approach to prove a conjecture of Fomin and Zelevinsky, that for cluster
algebras with no coefficients associated to quivers with no oriented
cycles, a seed is determined by its cluster.
We also obtain an interpretation of the monomial in the denominator of a
non-polynomial cluster variable in terms of the composition factors of an
indecomposable exceptional module over an associated hereditary algebra.
\end{abstract}

\maketitle

\section*{Introduction}

Cluster algebras were introduced by Fomin-Zelevinsky in \cite{fz1},
and have become increasingly important in many different parts of
algebra and beyond. 
We only consider the case involving only skew-symmetric matrices
and ``no coefficients''.
Essential
ingredients in the definition (in this case) are
certain $n$-element subsets $\underline{u}= \{ u_1, \dots , u_n\}$ of
the function field $\mathbb{Q}(x_1, \dots , x_n)$ called clusters, and pairs
$(\ul{u},Q')$ called seeds, where $\underline{u}$ is a cluster and
$Q'$ a finite quiver with no loops or 2-cycles (see Section 1 for
definitions). Here $Q'$ corresponds to a skew-symmetric integer matrix.

In this paper we only consider the acyclic case, 
i.e. cluster algebras for which there exists a seed $(\ul{u},Q')$, where $Q'$ 
is a quiver without oriented cycles. When we later refer to the {\em acyclic case}, 
the additional restrictions mentioned above are also assumed.

The aim of this paper is to prove that under our assumptions a seed
$(\ul{u}, Q')$ is determined by the cluster $\underline{u}$, as
conjectured by Fomin-Zelevinsky, who proved it for cluster algebras of
finite type \cite{fz2}.

Our method uses the theory of cluster categories and cluster-tilted
algebras from \cite{bmrrt}, \cite{bmr1} and \cite{bmr2}. These were
introduced to model some of the ingredients in the theory of cluster
algebras in a module theoretic/categorical way, inspired by \cite{mrz},
with the hope of
obtaining feedback to the theory of cluster algebras. The proof of 
our result
gives another illustration of this principle.

Suppose that $x$ is a non-polynomial cluster variable written in reduced form. Then 
the denominator is known to be a monomial \cite{fz1}.
An essential ingredient in our proof is to introduce the notion 
of {\em positivity condition}, which gives an easy criterion for an expression to be 
in reduced form.
Along the way we obtain an interpretation of the monomial in the denominator 
in terms of the composition factors of an indecomposable exceptional module
over an associated hereditary algebra, namely the path algebra of $Q$, where
$Q$ is a quiver with no oriented cycles appearing in a seed.  
Each indecomposable exceptional module corresponds to the 
denominator of a cluster variable this way.
This extends previous results
for finite type. The proof is also different from the earlier proofs in the
case of finite type. For finite type, we also obtain a new proof of the
existence of a bijection between cluster variables and indecomposable objects in the
corresponding cluster category inducing a bijection between clusters
and tilting objects, as first proved in \cite{bmrrt}.
In \cite{bmrrt} it was also conjectured that beyond
finite type, there should be 
a bijection between the
exceptional objects of the cluster category and the cluster variables.

That a seed is determined by its cluster in the 
acyclic case is also easily deduced from the proof of Theorem 4 in \cite{ck2}. 
(The result is stated in \cite{ck2} as a 
consequence of Theorem 4, which depends on this paper.)
While both approaches use heavily results from \cite{bmrrt} and \cite{bmr2},
there is an essential difference. We define a surjective map $\alpha$
from cluster variables to the indecomposable exceptional objects in the cluster
category, using our elementary notion of positivity. In \cite{ck2} a
surjective map $\beta$ is defined in the opposite direction by extending a beautiful formula
in \cite{cc} beyond finite type. Our approach gives the extra result on
denominators. However, combining the two 
surjectivity results gives directly stronger results, see Appendix.

Throughout this paper, $k$ denotes an algebraically closed field.

\section{background}
In this section we give some relevant background material for our
work.
\subsection{Cluster algebras}
We recall the main ingredients in the definition of cluster algebras
from \cite{fz1}, without dealing with the most general setting. We denote by
$\mathbb{Q}(x_1, \dots ,x_n)$ the function field in $n$ variables
over the rational numbers $\mathbb{Q}$, and by $Q$ a finite
quiver with vertices $1, \dots , n$, and with no oriented cycles. Then
$\underline{x}= \{x_1, \dots , x_n\}$ is a cluster, and we consider
the pair $(\underline{x},Q)$ to be an {\em initial seed}. For each
$i=1,\dots, n $, we denote the {\em mutation} of a quiver $Q$ by $\mu_i (Q)=Q'$.
We also denote the mutation of a seed 
$(\underline{x}, Q)$ by 
$\mu_i (\underline{x}, Q) = (\ul{x'} , Q')$. Here $x_i$ is
replaced by another element $x_i^{\ast}$, to get the cluster
$\underline{x}'$, where $x_i^{\ast}$ is defined by a
formula of the form $x_i \cdot x_i^{\ast} = m_1 + m_2$, with
$m_1$ and $m_2$ monomials depending on $\ul{x}$ and $Q$.
The precise definitions of $Q'$ and $x_i^{\ast}$ are given in \cite{fz1},
see also \cite{bmr2}.
Continuing this way we get a collection of
seeds. The $n$-element subsets of the seeds are the {\em clusters}, and the
elements in the clusters are the {\em cluster variables}. The associated
{\em cluster algebra} is the subalgebra $\A(Q)$ of
$\mathbb{Q}(x_1, \dots , x_n)$ over $\mathbb{Q}$ generated by the
cluster variables. We can start with any of the seeds as the initial
seed. The cluster algebras we have defined here are said to be {\em acyclic},
meaning that one of the quivers occurring in some seed, in this case $Q$ itself, has no
oriented cycles \cite{fz3}.

\subsection{Denominators of cluster variables}
When expressing a cluster variable in terms of the initial cluster
$\underline{x}$ as a reduced fraction 
$$f/m = f(x_1, \dots,x_n)/m(x_1, \dots ,x_n)$$ 
one has the remarkable fact that $m(x_1,\dots ,x_n)$ is a monomial. 
This is called the Laurent phenomenon
\cite{fz1}. For finite type, the quiver $Q$ in the initial seed
can be assumed to be a Dynkin quiver \cite{fz2}, and it is known that 
when a cluster variable different from $x_1, \dots, x_n$ is expressed in
terms of the cluster variables in the initial seed $(\ul{x}, Q)$ the denominator
determines the composition factors of an
indecomposable $kQ$-module. This was shown in \cite{fz2} for a Dynkin
quiver with alternating orientation, and extended to an arbitrary
Dynkin quiver in \cite{ccs2} (independently \cite{rt}). A
different proof was given in \cite{ck1} based upon \cite{cc}.

\subsection{Cluster categories}
For a finite dimensional hereditary $k$-algebra $H$, the
cluster category $\C_H= D^b(H)/\tau^{-1}[1]$ was
introduced in \cite{bmrrt} in order to model essential ingredients in
the theory of cluster algebras (see \cite{ccs1} for a category associated with 
the Dynkin diagram $A_n$, shown to be equivalent to the corresponding cluster category).
Here $D^b(H)$ denotes the
bounded derived category of $H$, with suspension functor $[1]$, 
and $\tau$ denotes the AR-translation
on $D^b(H)$, induced by the corresponding translation for
$H$-modules. 
The category $\C_{H}$ is known to be triangulated
\cite{k}. 
Any indecomposable object in $\C_H$ is represented 
by an object in $D^b(H)$ in the fundamental domain, that is, an 
indecomposable $H$-module or an object of the form $P[1]$ with $P$ 
an indecomposable projective $H$-module. 
We often identify an object in $\C_H$ with its representative in the
fundamental domain in $D^b(H)$.

Cluster-tilting theory was investigated in \cite{bmrrt},
and the {\em (cluster-)tilting objects} $T$ are the basic objects $T$ with
no self-extensions such that $T$ is maximal with this property, that
is, if $\Ext^1_{\C_{H}}(T\oplus X, T\oplus X)=0$,
then $X$ is in $\add T$. Any tilting object is induced by a tilting module over
some hereditary algebra derived equivalent to $H$. When $T=\bar{T} \oplus
M$ is a basic tilting object with $M$ indecomposable, there is a
unique indecomposable object $M^{\ast}\not \simeq M$ such that $T' =\bar{T}
\oplus M^{\ast}$ is a tilting object. Moreover, $M$ and $M^{\ast}$ are
related via {\em exchange-triangles} $M^{\ast} \xrightarrow{g} B\xrightarrow{f} M \to$
and $M \xrightarrow{t} B' \xrightarrow{s} M^{\ast} \to$ in $\C_H$ where $B$ and $B'$
are in $\add \bar{T}$, and $f$ and $s$ are minimal right $\add
\bar{T}$-approximations and $g,t$ are minimal left $\add
\bar{T}$-approximations. Here $(M,M^{\ast})$ is called an {\em exchange pair}.
The indecomposable objects occurring as direct summands of the 
tilting objects in $\C_{H}$ are exactly the indecomposable exceptional
objects of $\C_{H}$ (an object $X$ of $\C_{H}$ is said to be exceptional
if $\Ext^1_{\C_{H}}(X,X) = 0$).

\subsection{Cluster-tilted algebras}
The algebras $\End_{\C_{H}}(T)^{\op}$, where $T$ is a 
tilting object in $\C_{H}$ are called
{\em cluster-tilted algebras} and a study of
these algebras was initiated in \cite{bmr1}. In the acyclic
case it is known that if $Q$ is a quiver with no oriented cycles which
is part of a seed and $H=kQ$, then the quivers occurring in the
seeds of the cluster algebra $\A(Q)$ coincide with the
quivers of the cluster-tilted algebras associated with
$\C_{H}$ \cite{bmr2}. Also the multiplication rule $x_i \cdot
x_i^{\ast} = m_1 +m_2$ appearing in the definition of a cluster
algebra \cite{fz1} is interpreted in \cite{bmr2} in terms of the 
exchange triangles in the cluster category. 
The multiplication rule is independently interpreted in the
case of finite type in \cite{ccs2}.

\section{Discussion of the problem /motivation}\label{motiv}

We recall some results from \cite{bmr2} which are an essential starting
point, and which give motivation for our later considerations. 

Start with an initial seed $(\underline{x}, Q)$, where
$\underline{x}=\{x_1,\dots,x_n\}$ and $Q$ is a finite quiver without
oriented cycles with vertices $1,\dots,n$. Let $H=kQ$, and consider
the associated {\em tilting seed} $(H[1], Q)$. A tilting seed is a pair
$(T,\Gamma)$ where $T$ is a tilting object in $\C_H$ and $\Gamma$ is the
quiver of $\End_{\C_{H}}(T)^{\op}$. We can associate with the
cluster variable $x_i$ the indecomposable exceptional object $P_i[1]$,
where $P_i$ is the indecomposable projective $H$-module associated
with the vertex $i$. Applying the mutation $\mu_i$ to
$(\underline{x},Q)$ we get $(\underline{x}',Q')$, with
$\underline{x}'$ containing a new cluster variable $x_i^{\ast}$, 
replacing $x_i$. Define
$\delta_i$ of $(H[1],Q)$ to be $(T',Q'')$ where $T'$ is obtained from
$H[1]$ by replacing $M=P_i[1]$ with a non-isomorphic indecomposable
object $M^{\ast}$ to get a new tilting object $T'$. Here $Q''$ is by
definition the quiver of $\End_{\C_{H}}(T')^{\op}$. We then
associate $M^{\ast}$ with $x_i^{\ast}$. Thus the elements of the
new cluster $\underline{x}'$ are associated with the summands of the tilting
object $T'$. We also know that the quivers $Q'$ and $Q''$ are
isomorphic \cite{bmr2}. Continuing with a sequence of corresponding
mutations for the seeds and for the tilting seeds, we associate
in this way an indecomposable exceptional object in $\C_{H}$
to each cluster variable occurring for the chosen sequence of mutations.
Note that, a priori, the exceptional object associated to a
cluster variable is not uniquely determined by the cluster variable itself
--- it depends on the sequence of mutations chosen.
Clusters are taken to tilting objects, and at each stage the corresponding
quivers coincide.

Recall that we want to prove that a seed is determined by its
cluster. The corresponding result for tilting seeds is true by
definition, since the corresponding quiver is given by the
endomorphism algebra of the tilting object. So it should be
useful to get a closer relationship between the seeds and the tilting
seeds. 
In particular, we show in this paper that in the above approach, the exceptional 
object associated to a cluster variable does not depend on the sequence 
of mutations used to reach it. Then a cluster $\ul{x}'$ determines a tilting 
object $T'$ uniquely and hence the quiver $Q'$ of the tilting seed $(T',Q')$ 
and therefore the quiver of the seed which has cluster $\ul{x}'$.
Note that if there is a map from
the set of cluster variables in $\A(Q)$ to the set
of indecomposable exceptional objects in $\C_H$, taking
$x_i$ to $P_i[1]$, and sending clusters to tilting objects,
then it is uniquely determined.

So the strategy is to get enough information on the shape of each
cluster variable in order to associate to it an indecomposable
exceptional object in $\C_{H}$, in a way compatible with the
above locally defined map. We can then use this to improve the results
from \cite{bmr2} and prove the desired conjecture. This will be based
upon proving two theorems, stated below, which should be of independent
interest. Before we state them, we introduce some more terminology.

Let $H=kQ$ and let $S_1,\dots,S_n$ denote the simple $H$-modules
associated with the vertices $1,\dots,n$. For $M$ in $\mod H$ denote
by $[M]$ the corresponding element in the Grothendieck group
$K_0(H)$. For a monomial $m=m(x_1,\dots,x_n)$, denote by $\beta(m)$
the corresponding element in $K_0(H)$, induced by
$\beta(x_i^{n_i})=n_i[S_i]$. For a cluster variable $u$, which is not a polynomial,
we write
$u =f/m$ in reduced form, where $f$ is a polynomial and $m$ is a non-constant monomial with
constant term 1. If $u$ is a polynomial and $x_i \mid u$ 
for some $x_i $ in $\underline{x}$, we write $u =x_i f$. We say
that a polynomial $f=f(x_1,\dots,x_n)$ satisfies the {\em positivity
condition} if $f(e_i)>0$ where $e_i=(1,\dots,1,0,1,\dots,1)$
(with 0 in the i'th position) for $i=1,\dots,n$.
The following is a crucial observation.

\begin{lem}
Let $u=f/m$, where $f$ is a polynomial satisfying the positivity 
condition and $m$ is a non-constant monomial. Then $u=f/m$ is in reduced form.
\end{lem}

We say that a cluster variable expressed in terms of the initial seed 
$(\underline{x}, Q)$, where $Q$
is a quiver without oriented cycles, is of type ($\star$) if 
\begin{itemize}
\item[-]{It is of the form $f/m$ or $fx_i$
where $m$ is a non-constant monomial and $f$ is a polynomial
satisfying the positivity condition.}
\item[-]{In case it is equal to $f/m$, 
there is an exceptional $H$-module $M$, such that $\beta(m)= [M]$.}
\end{itemize}
Recall that an indecomposable exceptional $H$-module $M$ is 
uniquely determined by $[M]$ (see \cite{ker}). Thus, the $M$
appearing in the above definition is uniquely defined.

Our first main result is the following.

\begin{thm}\label{thm1}
The cluster variables for the cluster algebra $\A(Q)$, are all of type ($\star$).
For each non-polynomial cluster variable $f/m$, there is a unique 
indecomposable exceptional $H$-module $M$ such that $\beta(m) =[M]$, 
and each indecomposable exceptional module corresponds to the 
denominator of some cluster variable.
\end{thm}

We define a map $\alpha$ from cluster variables of type ($\star$) to indecomposable 
objects in $\C_H$, with $H= kQ$, as follows. 
\begin{enumerate}
  \item[(i)]
    $\alpha(f/m)=M$, where $M$ is the indecomposable 
exceptional object with $\beta(m) =[M]$.
  \item[(ii)]
    $\alpha(fx_i)=P_i[1]$, where $P_i$ is the projective cover of $S_i$.
\end{enumerate}
Let $c_M = \beta^{-1}[M]$, that is, the exponents of $c_M$ give the
composition factors of $M$.

There is an induced map $\bar{\alpha}$
defined on the clusters consisting of cluster variables of type ($\star$)
to the set of tilting objects, and an induced map
$\widetilde{\alpha}$ from seeds with the same condition on the involved clusters,
to the set of tilting seeds.

The map $\alpha$ is used in the proof of Theorem \ref{thm1}, and
it is clear that when this proof is finished, the map $\alpha$ is
defined on the set of cluster variables, and hence
$\bar{\alpha}$ and $\widetilde{\alpha}$ are defined on the 
sets of clusters and seeds, respectively.

Note that if $X$ is an indecomposable $H$-module, then by \cite{bmrrt} we have
$$\Ext^1_{\C_{H}}(X,X) \simeq \Ext^1_{
H}(X,X)\oplus D \Ext^1_{H}(X,X),$$ so that $X$ is an exceptional
$H$-module if and only if it is an exceptional object in
$\C_{H}$. Also $P_i[1]$ is an indecomposable exceptional object in
$\C_{H}$.

The second theorem is the following.

\begin{thm}\label{thm2}
Let $Q$ be a finite quiver with no oriented cycles, and let $H=kQ$.
The map $$\alpha \colon \{\text{cluster variables of $\mathcal{A}(Q)$} \} \to
  \{\text{indecomposable exceptional objects in }\C_{H}\}$$ is
  surjective, it induces a surjective map $\bar{\alpha} \colon \{
  \text{clusters}\} \to \{ \text{tilting objects}\}$, and
a surjective map $\widetilde{\alpha} \colon \{
  \text{seeds}\} \to \{ \text{tilting seeds}\}$, preserving quivers.
\end{thm}

Those two results will actually be proven together, where the crucial
induction step is treated in the next section.

\section{The induction step}
In this section we give the induction step needed for the proof of
Theorems \ref{thm1} and \ref{thm2}. As before $(\underline{x},Q)$, with
$\underline{x}=\{x_1,\dots,x_n\}$, denotes the initial seed, with $Q$
a finite quiver with no oriented cycles and vertices $ 1, \dots, n$.

\begin{prop}\label{ind}
Let $(\ul{x}',Q')$ be a seed, with $\ul{x}' = \{x_1', \dots, x_n' \}$,
and choose some $x_i'$. Assume that all $x_j'$ are of type ($\star$), and that
$\alpha(x_1') \oplus \cdots \oplus \alpha(x_n') = T' = M \oplus \bar{T}$ is a tilting object 
in $\C_H$, where
$M = \alpha(x_i')$, and $Q'$ is the quiver of $\End_{\C_H}(T')^{\op}$. For an indecomposable
direct summand $A$ of $T'$, denote the corresponding cluster variable in $\ul{x}'$ by $x_A$.
Then we have the following.
\begin{itemize}
\item[(a)]{Performing the mutation $\mu_i$ on $(\ul{x}', Q')$, the new cluster variable 
$(x_i')^{\ast} = (x_M)^{\ast}$ is also of type ($\star$), and $\alpha((x_M)^{\ast})= M^{\ast}$, where
$M^{\ast}$ is the unique indecomposable object in $\C_H$ with $M^{\ast} \not \simeq M$ and 
$\bar{T} \oplus M^{\ast}$ a tilting object.}
\item[(b)]{If $(\ul{x}'', Q'')$ is the new seed, where $\ul{x}'' 
= \{x_1', \dots, (x_i')^{\ast}, \dots, x_n' \}$,
then $$\alpha(x_1') \oplus \cdots \oplus \alpha((x_i')^{\ast}) \oplus \cdots \oplus \alpha(x_n') 
= \bar{T} \oplus M^{\ast},$$
so that the cluster $\ul{x}''$ is taken to the tilting object $T'' = \bar{T} \oplus M^{\ast}$, 
and $Q''$ is the quiver
of $\End_{\C_H}(T'')^{\op}$, so that the seed $(\ul{x}'', Q'')$ is taken to the tilting seed
$(T'', Q_{T''})$, where $Q'' = Q_{T''}$.
}
\end{itemize}
\end{prop}

\begin{proof}
Let $B \to M$ and $M \to B'$ be minimal right and left $\add \bar{T}$-approximations
of $M$, respectively. If $A$ is an 
indecomposable direct summand of $M \oplus B \oplus B'$, then it is
an indecomposable direct summand of $T'$. Hence we have by assumption that the cluster variable $x_A$
is of type ($\star$), that is, 
$x_A = f_A/m_A$, with $\beta(m_A) = A$ if $A$, chosen in the fundamental domain, is an $H$-module,
and $x_A = f_A \cdot x_i = f_A \cdot c_{S_i}$ if $A = P_i[1]$
with $S_i$ the simple module corresponding to the vertex $i$, and
$P_i$ the projective cover of $S_i$. In both cases 
$f_A$ satisfies the positivity condition.

Consider the approximation triangles
$$M^{\ast} \to B \to M \to \text{  and  } M \to B' \to M^{\ast} \to,$$
where $B= \oplus_i B_i$, and $B= \oplus B_j'$, with the $B_i$ and $B'_j$ indecomposable.

We have by \cite{bmr2} that $$x_M \cdot (x_M)^{\ast} = x_B + x_{B'},$$
where
$x_B = \prod x_{B_i}$ and $x_{B'} = \prod x_{B'_j}$.
For the further analysis we consider two different cases.

\vspace{0.7cm}

\noindent CASE I.

\noindent Assume that we have an exact sequence $0 \to M^{\ast} \to B \to M \to 0$ of $H$-modules.
Since $M^{\ast}$ is not injective, $\tau^{-1}M^{\ast}$ is an $H$-module.
By \cite{bmrrt}, we have that 
$D \Ext^1_{\C_H}(M,M^{\ast}) \simeq \Hom_{\C_H}(M^{\ast}, \tau M) \simeq \Hom_{\C_H}(\tau^{-1}M^{\ast},M)$
is $1$-dimensional over $k$.
Since $\Hom_{D^b(H)}(\tau^{-1} M^{\ast}, \tau^{-1}M[1]) = \Ext^1_H(M^{\ast},M) = 0$,
because $\Ext^1_H(M,M^{\ast}) \neq 0$, there is a non-zero module map $h \colon \tau^{-1} M^{\ast} \to M$.
Consider the corresponding exact sequence
$$0 \to K_h \to \tau^{-1} M^{\ast} \overset{h}{\to} M \to L_h \to 0.$$
This gives rise to the triangle 
 $$M \to L_h \oplus K_h[1] \to M^{\ast} \to $$
in $D^b(H)$, whose image in $\C_H$ is (isomorphic to) the exchange triangle
$$M \to B' \to M^{\ast}\to.$$
Write $K_h = P \oplus K_h'$, where $P$ is projective and 
$K_h'$ has no non-zero projective direct summands.
Then, in the cluster category, we have
$B' \simeq L_h \oplus K_h'[1]  \simeq L_h \oplus \tau K_h' \oplus P[1]$.
Here, $ L_h \oplus \tau K_h' \oplus P[1]$ is in the fundamental domain in $D^b(H)$.
Since we have a monomorphism $K_h \to \tau^{-1}M^{\ast}$, we also have a monomorphism
$\tau K_h' \to M^{\ast}$. Letting $A = L_h \oplus \tau K_h'$, we see that
 $[B] = [M] + [M^{\ast}] = [A] + [C]$, for some $H$-module $C$, using that 
$L_h$ is a factor of $M$, and that $\tau K_h'$ is a submodule of $M^{\ast}$.

For a direct summand $Y = Y_1 \oplus \cdots \oplus Y_r$ of $T'$, where the $Y_j$ are indecomposable,
we write $f_Y = f_{Y_1} \cdots f_{Y_r}$, $m_Y = m_{Y_1} \cdots m_{Y_r}$ and
$x_Y= x_{Y_1} \cdots x_{Y_r}$. Then we have $x_{B'} = x_A \cdot x_{P[1]} = (f_A/m_A) \cdot c_S \cdot f_{P[1]}$,
where $S= P/\r P$. Here $f_A$ satisfies the positivity condition, since $f_{A'}$ satisfies the positivity 
condition for each indecomposable direct summand $A'$ of $A$.
We have $m_B = m_A \cdot m'$, where $m'$ is a monomial, using that we know that
$m_B=c_B$ and $m_A=c_A$.

By the interpretation of the multiplication $x_M \cdot (x_M)^{\ast}$ 
in the cluster category from \cite{bmr2}, we have that
\begin{equation}\label{one}
x_M \cdot (x_M)^{\ast} = x_B + x_{B'} = \frac{f_B}{m_B} + \frac{f_A}{m_A} \cdot c_S \cdot f_{P[1]},
\end{equation} 
so that $$(x_M)^{\ast} = \frac{(f_B +f_A \cdot c_S \cdot f_{P[1]} \cdot m') \cdot m_M}{m_B \cdot f_M}.$$
It follows from the exact sequence $0 \to M^{\ast} \to B \to M \to 0$ that
$m_B/m_M = c_{M^{\ast}}$, the monomial corresponding to the composition factors of $M^{\ast}$.
We have
$$(x_M)^{\ast} = \frac{(f_B +f_A \cdot c_S \cdot f_{P[1]} \cdot m') \cdot \frac{1}{f_M}}{c_{M^{\ast}}},$$
Note that as $f_M$ is assumed to satisfy the positivity condition, it can
have no monomial factor, and it follows from the Laurent phenomenon that
the numerator is a polynomial.
We want to show the desired positivity condition for the numerator,
which in particular will imply that the above expression for $(x_M)^{\ast}$ is in reduced form:
For $e_i = (1, \cdots,1,0,1,\cdots,1)$ (the $0$ in the $i$'th position) we have by assumption, or
directly, that $f_B(e_i) > 0$, $f_A(e_i)>0$, $f_{P[1]}(e_i) > 0$, $m'(e_i) \geq 0$, 
$c_S(e_i) \geq 0$ and $f_M(e_i)>0$,
so overall we get a positive number. The expression is thus in its reduced
form, and hence we get
$m_{M^{\ast}} = c_{M^{\ast}}$ and consequently $\alpha((x_M)^{\ast}) = M^{\ast}$, as desired.

If $M$ and $M^{\ast}$ are still assumed to be $H$-modules, and there is no non-split exact sequence
$0 \to M^{\ast} \to B \to M \to 0$ of $H$-modules, then we have a non-split exact sequence 
$0 \to M \to B' \to M^{\ast} \to 0$ of $H$-modules, by \cite{bmrrt}. This amounts to,
using the previous notation, assuming our claims are proved for $B,B'$ and $M^{\ast}$, and then 
showing the desired property for $M$. So the equation (\ref{one}) 
is replaced by 
\begin{equation}\label{onep}
x_{M^{\ast}} \cdot (x_{M^{\ast}})^{\ast} = \frac{f_B}{m_B} + \frac{f_A}{m_A} \cdot c_S \cdot f_{P[1]}.
\end{equation} 
We then get 
$$(x_{M^{\ast}})^{\ast} = \frac{(f_B +f_A \cdot c_S \cdot f_{P[1]} \cdot m') \cdot \frac{1}{f_{M^{\ast}}}}{c_M},$$
and we get similarly that $c_M = m_M$ and $\alpha((x_{M^{\ast}})^{\ast})= M$.

\vspace{0.7cm}
\noindent CASE II.

\noindent Assume now that $M$ is not induced by a module, hence $M = P[1]$ for an indecomposable projective $H$-module $P$.
Since $\Ext^1_{\C_H}(M,M^{\ast}) \neq 0$, we have that $M^{\ast}$ cannot be of the form $Q[1]$,
with $Q$ projective, and is hence an $H$-module.

Consider the diagram
$$
\xy
\xymatrix{
& & K_g \oplus I'[-1] \ar[d]& & \\
B_1[-1] \oplus P'' \ar[r] & P \ar[r] & M^{\ast} \ar[r] \ar[d]_{g} & B_1 \oplus P''[1] \ar[r]^{f} & P[1] \\
& & I \ar[d] & & \\
& & I' \oplus K_g[1] & & 
}
\endxy
$$
where $f$ is a minimal right $\add \bar{T}$-approximation of $M = P[1]$, and $B = B_1 \oplus P''[1]$,
where $B_1$ is an $H$-module and $P''$ a projective $H$-module, $I=\tau(P[1])$
is the indecomposable injective $H$-module
with socle isomorphic to the simple module $S= P/ \r P$, 
and $g \colon M^{\ast} \to I$ a non-zero map in the cluster category, and hence in the module
category, since $I$ is injective.
We have that $\Hom_{\C}(M^{\ast},I) \simeq \Hom_{H}(M^{\ast},I)$ is 1-dimensional, since
$\Hom_{\C}(M^{\ast},I) = \Hom_{\C}(M^{\ast}, \tau P[1]) = \Hom_{\C}(M^{\ast},\tau M) = \Hom_{\C}(M^{\ast},M[1])$
and $ \Hom_{\C}(M^{\ast},M[1])$ is 1-dimensional. 
We have in addition that $I[-1] \simeq P[1]$ in the cluster category.
Hence it is clear that the image of the triangle induced from the map $g$ is indeed isomorphic
to the exchange triangle $M \to B' \to M^{\ast} \to$.
This means that we have the minimal right
$\add \bar{T}$-approximation $B'= K_g \oplus I'[-1] \to M^{\ast}$ of $M^{\ast}$, where $I' = \Coker g$ and
$K_g = \Ker g$. Here $I'$ is an injective $H$-module, being a factor of the injective $H$-module $I$.
Denote its socle by $S'$, and let $P'$ be the projective cover of $S'$.
We have exact sequences $0 \to P'' \to P \to M^* \to B_1 \to 0$
and $0 \to K_g \to M^{\ast} \to I \to I' \to 0$. Consider the exact commutative diagram

$$
\xy
\xymatrix{
& 0 \ar[d] & 0 \ar[d] & & \\
0 \ar[r] & \r P / P'' \ar[d] \ar[r] & K_g \ar[d] \ar[r] & C \ar[d]_{u} \ar[r] & 0\\
0 \ar[r] & P / P'' \ar[d] \ar[r] & M^{\ast} \ar[r] & B_1 \ar[r] & 0 \\
& S \ar[d] & & & \\
& 0 & & & 
}
\endxy
$$
We use here that $\r P/P'' \subseteq K_g$, since the composition $P \to M^{\ast} \to I$
has image in the simple socle $S$ of $I$. Since $\dim_k \Hom(M^{\ast},I) = 1$, we have that $S$
occurs exactly once as a composition factor of $M^{\ast}$. It also occurs exactly once as a
composition factor of $I$ as there are no loops in the quiver of $H$. 
The exact sequence $0 \to K_g \to M^{\ast} \to I$ then shows that
$S$ is not a composition 
factor of $K_g$, and hence not of $C$ or of $\Ker u$. By the snake lemma applied to 
the above diagram we have an exact sequence $0 \to \Ker u \to S$, which
shows that $u \colon C \to B_1$ is a monomorphism.

We have by assumption that the cluster variables in $\ul{x'}$ corresponding to
indecomposable direct summands of 
$B,B'$ and $M$ are of type ($\star$).
Hence we have
$x_{P[1]} = f_{P[1]} \cdot x_{S}$ and $x_{P''[1]} = f_{P''[1]} \cdot x_{S''}$
and $x_{P'[1]} = f_{P'[1]} \cdot x_{S'}$,
with $f_{P[1]}$, $f_{P'[1]}$ and $f_{P''[1]}$ satisfying the positivity condition, and
$S''=P''/\r P''$.
We also have $m_{K_g} = c_{K_g}$ and $m_{B_1} = c_{B_1}$. This gives
\begin{align*}
x_M \cdot (x_M)^{\ast} & = 
x_B + x_{B'}  \\
& = x_{B_1} \cdot x_{P''[1]} + x_{K_g} \cdot x_{P'[1]} \\
& = \frac{f_{B_1}}{c_{B_1}} \cdot c_{S''} \cdot f_{P''[1]} + \frac{f_{K_g}}{c_{K_g}} \cdot c_{S'} \cdot f_{P'[1]}.
\end{align*}
By the above exact sequences we have $c_{B_1} = c_{B_1/C} \cdot c_C$ and
$c_{K_g} = c_{\r P / P''} \cdot c_C$ and $[M^{\ast}] = [P/P''] + [B_1] = [\r P/P''] + [S] 
+[B_1]$, so that $c_{M^{\ast}}/c_S = c_{\r P/P''} \cdot c_{B_1/C} \cdot c_C$.
This gives, using $x_M = x_{P[1]} = f_{P[1]} \cdot x_S$, that
$$(x_M)^{\ast} = \frac{(c_{\r P/P''} \cdot f_{B_1} \cdot f_{P''[1]} \cdot c_{S''} + 
c_{B_1/C} \cdot f_{K_g} \cdot f_{P'[1]} \cdot c_{S'}) \cdot \frac{1}{f_M}}{c_{M^{\ast}}}.$$

We want to show that $\r P / P'' \oplus S''$ and $B_1/C \oplus S'$ 
have no common composition factors, so that no $x_i$
divides $c_{\r P/ P''} \cdot c_{S''}$ and $c_{B_1/C} \cdot c_{S'}$.
If $S_1$ is a simple $H$-module which is a composition factor of $\r P/P'' \oplus S''$,
there must be a path in the quiver of $H$ from the vertex corresponding to $S$
to the vertex corresponding to $S_1$. 
We have $M^{\ast}/K_g \subseteq I$, and $I'$ is a factor of $I$.
If $S_1$ is a composition factor of $B_1/C \oplus S'$, and hence of
$M^{\ast}/K_g \oplus S'$, we would have a path from the vertex corresponding 
to $S_1$ to the vertex corresponding to $S$.
Since the quiver of $H$ has no oriented cycles, we conclude that 
no $x_i$ divides both $c_{\r P/P''} \cdot c_{S''}$ and 
$c_{B_1/C} \cdot c_{S'}$. 
Since $f_{B_1}$, $f_{K_g}$, $f_M$, $f_{P''[1]}$ and
$f_{P'[1]}$ satisfy the positivity condition, it follows that the numerator
in the above expression for $(x_M)^{\ast}$ also does.
Since the denominator gives the composition factors for the indecomposable exceptional 
$H$-module $M^{\ast}$, the cluster variable
$(x_M)^{\ast}$ is of the desired form, and we have 
$\alpha((x_M)^{\ast}) = M^{\ast}$.

The case where $M^{\ast} = P[1]$ and $M$ is an $H$-module follows by symmetry as in Case I.
In view of the discussions in Section 2, this finishes the proof.
We have then proved (a). The first part of (b) then follows directly,
and it follows from \cite{bmr2}, as discussed in Section \ref{motiv},
that $(\ul{x}'', Q'')$ and $\tilde{\alpha}(\ul{x}'', Q'')$ have the same quiver.
\end{proof}

\section{Proof of main result}
In this section we apply the induction step to give the proof of Theorems 
\ref{thm1} and \ref{thm2}, and we give our application to seeds being
determined by clusters.

\begin{proof}[Proof of Theorems \ref{thm1} and \ref{thm2}]
Let $(\ul{x},Q)$, with $\ul{x} = \{x_1, \dots, x_n \}$, be an initial seed, where
$Q$ has no oriented cycles. Then $x_1, \dots , x_n$ are clearly all of the desired form, and $\alpha(x_i)$
is defined to be $P_i[1]$. 
Then $\alpha(x_1) \oplus \cdots \oplus \alpha(x_n) = P_1[1] \oplus \cdots \oplus P_n[1]$ is a tilting object 
in $\C_H$, with $H= kQ$, so the seed $(\ul{x}, Q)$ corresponds to the tilting seed $(H[1], Q)$,
which is the start of our induction.

Consider the following property ($\widetilde{\star}$) of a seed $(\ul{u}, \G)$, where
$\ul{u} = \{u_1, \dots, u_n \}$.

\noindent {\bf Property} ($\widetilde{\star}$): The cluster variables $u_1, \dots , u_n$
are of type $(\star)$, so that $\alpha(u_i)$ is defined for $i= 1, \dots, n$,
the induced map $\bar{\alpha}$
takes the cluster $\ul{u}$ to the tilting object $T = \alpha(u_1) \oplus \cdots \oplus \alpha(u_n)$,
and the induced map $\widetilde{\alpha}$ takes the seed $(\ul{u}, \G)$ to the tilting seed
$(T,Q_T)$, where $\G = Q_T$.

By the above remarks, the initial seed $(\ul{x},Q)$ satisfies ($\widetilde{\star}$).
Let $(\ul{x}'',Q'')$ be a seed where $t>0$ is the number of mutations in
a minimal sequence of mutations from $(\ul{x},Q)$ to $(\ul{x}'',Q'')$.
Assume that all seeds reached by a sequence of mutations in a minimal sequence of mutations
of length less that $t$ have property ($\widetilde{\star}$). Then choose $(\ul{x}',Q')$, such
that $(\ul{x}'',Q'') = \mu_i (\ul{x}',Q')$ for some $i$, and $(\ul{x}', Q')$
is reached by a minimal sequence of mutations of length $t-1$.
Then it follows directly from Proposition \ref{ind}
that $(\ul{x}'',Q'')$ satisfies ($\widetilde{\star}$).

It then only remains to show that $\alpha$, $\bar{\alpha}$ and $\widetilde{\alpha}$ are surjective.
Let $T'_i$ be an exceptional object in $\C_H$. It is then a direct summand of a tilting object $T'$.
Since the cluster-tilting graph is connected (see \cite{bmrrt}), 
the corresponding tilting seed $(T', Q')$ is obtained from
$(H[1],Q)$ by a finite number of mutations of tilting seeds. We do the corresponding
sequence of mutations of seeds to the initial seed $(\ul{x},Q)$.
Then we are in the above situation, and it follows that all the maps $\alpha$, $\bar{\alpha}$ and
$\widetilde{\alpha}$ are surjective.
\end{proof}

The following was conjectured by Fomin and Zelevinsky \cite{fz2} (for
any cluster algebra).

\begin{thm}
For an acyclic cluster algebra (associated with a skew-symmetric matrix and with no coefficients) 
a seed is determined by its cluster.
\end{thm}

\begin{proof}
We can choose an initial cluster $(\ul{x},Q)$, where $Q$ is a quiver with no oriented cycles,
and express the cluster variables in terms of the variables $x_1, \dots, x_n$ in $\ul{x}$.
Let $\ul{x}'$ be a cluster, and choose a seed $(\ul{x}',Q')$. Consider the tilting seed
$\widetilde{\alpha}(\ul{x}',Q') = (\bar{\alpha}(\ul{x}'),Q')$. Since $\bar{\alpha}(\ul{x}')$
determines $Q'$, it follows that the cluster $\ul{x}'$ determines the seed $(\ul{x}',Q')$.
\end{proof}

\section{Concluding remarks}

Note that for finite type it is known \cite{fz2} that the number
of cluster variables is equal to the number of
almost positive roots, and hence to the number of indecomposable objects (which are all exceptional)
in the associated cluster category .
Hence knowing that $\alpha$ is surjective, we can conclude from this
that $\alpha$ is actually a bijection. As a consequence, this gives an alternative proof of the correspondence theorem, 
taking clusters to tilting objects, in \cite{bmrrt}.

If one is only interested in the shape of the denominator for cluster 
variables for finite type, one can get a more direct proof, still using 
our notion of positivity: define a map $\psi \colon \ind \C_H \to \{ \text{cluster variables} \}$
by sending $P_i[1] =I_i[-1]$ to $x_i$ as usual, and using the almost split triangles
corresponding to the exchange pairs $(I_i[-1],P_i)$ 
to get cluster variables of the desired form
corresponding to the indecomposable projectives. 
Then use 
almost split sequences to get cluster variables of the desired form
corresponding to all the other objects. Then $\psi$ will clearly be 1-1, and
hence bijective by the above mentioned result from \cite{fz2}.

When $H$ is an indecomposable hereditary $k$-algebra of infinite type, we can use
the same procedure to get a 1--1 map $\psi$ from the indecomposable preprojective 
$H$-modules to the cluster variables. This gives in our setting a short proof for the fact that
if the quiver of $H$ is not Dynkin, the associated cluster algebra is of infinite type, which is
part of \cite[Th. 1.5]{fz2}.

\section*{Appendix A: Clusters and tilting objects for acyclic cluster algebras}

\noindent By Aslak Bakke Buan, Philippe Caldero, Bernhard Keller, Robert J. Marsh, Idun Reiten, 
Gordana Todorov
\vspace{5mm}
 
Let $H$ be a finite dimensional hereditary algebra over an algebraically closed field 
$k$ with associated cluster category $\C_H$ and cluster algebra
$\A_H$. Let $\ind \exc \C_H$ denote the set of isomorphism classes of indecomposable
execptional objects in $\C_H$.
The purpose of this appendix is to point out that combining the surjective 
Caldero-Chapoton map $\beta$ from \cite{ck2}, with the surjective map $\alpha$
in this paper gives the following stronger result, confirming Conjecture 9.1 in
\cite{bmrrt}, as a direct consequence.
It gives a refinement of Theorem 4 in \cite{ck2}, which is proved using the surjective 
map $\beta$, together with the idea of positivity from this paper.

\begin{thmA}
The surjective maps $$\alpha \colon \{\text{cluster variables} \} \to 
\{\ind \exc \C_H \}$$ and 
$$\beta \colon 
\{\ind \exc \C_H \} \to  \{\text{cluster variables} \} $$
are inverse bijections, inducing inverse bijections between 
clusters and tilting objects and between seeds and tilting seeds. 
\end{thmA}

\begin{proof}
Consider the surjective map 
$$\alpha \beta \colon 
\{ \ind \exc \C_H\} \to \\ 
\{ \ind \exc  \C_H \}. $$
It is clear from
the definition of $\alpha$ and $\beta$, that $\alpha \beta$
is the identity on the indecomposable direct summands of $H[1]$. It follows from
Theorem \ref{thm2} and the part of Theorem 4 in \cite{ck2}
that does not depend on this paper, that $\alpha \beta$ takes tilting objects
to tilting objects and tilting seeds to tilting seeds. It then follows 
from the recurrence relations satisfied by
the maps that
$\alpha \beta$ is the identity, and hence $\alpha$ and $\beta$ are inverse bijections.
\end{proof}

The following improvement of Theorem \ref{thm1} follows directly.

\begin{cor1}
The maps $\alpha$ and $\beta$ give inverse 1--1 correspondences between the
non-polynomial cluster variables and the isomorphism classes 
of indecomposable $H$-modules, such that the denominator of a cluster
variable is given by the composition factors of the corresponding indecomposable
$H$-module.
\end{cor1}

Another direct consequence is the following generalization of a result
proved for finite type in \cite{fz2}.

\begin{cor2}
Let $\underline{u} = \{u_1, u_2, \dots , u_n \}$ be a cluster for $\A_H$.
Then for each $i= 1, \dots, n$ there is exactly one cluster variable
$u_i^{\ast} \neq u_i$ such that $$\{ u_1, \dots, u_{i-1}, u_i^{\ast}, u_{i+1}, \cdots, 
u_n \}$$ is a cluster.
\end{cor2}

\end{document}